\documentclass[11pt]{article}

\usepackage{amsmath, amssymb, amsthm, verbatim,enumerate,bbm,color} 

\usepackage{mathrsfs}

\usepackage[colorlinks,citecolor=blue,bookmarks=true]{hyperref}

\title{On Erd\H{o}s's Method for Bounding the Partition Function}

\author{Asaf Cohen Antonir\thanks{School of Mathematical Sciences, Tel Aviv University, Tel Aviv, 6997801, Israel. Email: asafc1$@$tauex.tau.ac.il} \and Asaf Shapira \thanks{
School of Mathematics, Tel Aviv University, Tel Aviv 69978, Israel.
Email: asafico$@$tau.ac.il. Supported in part by ISF Grant 1028/16, ERC Starting Grant 633509 and NSF-BSF Grant 2019679.}}

\date{\today}
\allowdisplaybreaks

\theoremstyle{plain}
\newtheorem{theorem}{Theorem}

\makeatletter
\def\moverlay{\mathpalette\mov@rlay}
\def\mov@rlay#1#2{\leavevmode\vtop{%
   \baselineskip\z@skip \lineskiplimit-\maxdimen
   \ialign{\hfil$\m@th#1##$\hfil\cr#2\crcr}}}
\newcommand{\charfusion}[3][\mathord]{
    #1{\ifx#1\mathop\vphantom{#2}\fi
        \mathpalette\mov@rlay{#2\cr#3}
      }
    \ifx#1\mathop\expandafter\displaylimits\fi}
\makeatother



\RequirePackage{color}\definecolor{RED}{rgb}{1,0,0}\definecolor{BLUE}{rgb}{0,0,1} 

\begin{document}
\date{}
\maketitle

\begin{abstract}
For fixed $m$ and $R\subseteq \{0,1,\ldots,m-1\}$, take $A$ to be the set of positive integers congruent modulo $m$ to one of the elements of $R$,
and let $p_A(n)$ be the number of ways to write $n$ as a sum of elements of $A$. Nathanson proved that $\log p_A(n) \leq (1+o(1)) \pi \sqrt{2n|R|/3m}$
using a variant of a remarkably simple method devised by Erd\H{o}s in order to bound the partition function.
In this short note we describe a simpler and shorter proof of Nathanson's bound.

\end{abstract}

\section{Introduction.}\label{sec:intro}
A partition of an integer $n$ is a sequence of positive integers $a_1 \leq a_2 \leq \cdots$ whose sum is $n$.
Let $p(n)$ denote the classical partition function of $n$, namely, the number of ways to write $n$ as a sum of positive integers.
The celebrated Hardy--Ramanujan formula \cite{HR1} (discovered independently by Uspensky \cite{Usp}) states that $p(n) \sim \frac{1}{4n\sqrt{3}} \exp(\pi\sqrt{2n/3})$.
Erd\H{o}s \cite{Erdos} later devised a remarkably simple proof of the slightly weaker upper bound
\begin{equation}\label{eqErdos}
\log p(n) \leq \pi\sqrt{2n/3}\;.
\end{equation}
Let $\mathbb {N}$ denote the set of positive integers, and suppose $S\subseteq \mathbb {N}$. We define $p_S(n)$ to be the number of partitions of $n$ with all summands in $S$.
For a fixed positive integer $m$ and $R\subseteq \{0,1,\ldots,m-1\}$, we take $A=A(m,R)$ to be the set of all positive integers $a$ with $a \;(\mathrm{mod}\; m) \in R$.
Nathanson \cite{Nat1999} used Erd\H{o}s's method for proving \eqref{eqErdos} to obtain\footnote{Nathanson \cite{Nat1999} also proves that $\log p_A(n) \geq (1-o(1)) \pi\sqrt{2n|R|/3m}\;$.}
\begin{equation}\label{eqNathanson}
    \log p_A(n) \leq (1+o(1)) \pi\sqrt{2n|R|/3m}\;.
\end{equation}
The argument in \cite{Nat1999} was more complicated than Erd\H{o}s's  due to the need to control various error parameters (but was
still simpler than the original proof of this result \cite {Meinardus}); see the remark at the end of the proof.

Our goal in this short note is to give a proof of \eqref{eqNathanson} which is as simple as Erd\H{o}s's proof of (\ref{eqErdos}).
The main trick is that, instead of directly bounding $p_A(n)$, we will instead bound $p_{A^+}(n)$, where given $m$ and $R$ as above, we take $A^{+}=A\setminus R$, that is, the set of all integers $a \geq m$ with $a \;(\mathrm{mod}\; m) \in R$.
Our main result here is the following generalization\footnote{Indeed, when $m=1$ and $R=\{0\}$, we have $p_{A^+}(n)=p(n)\;$.} of \eqref{eqErdos}.

\begin{theorem}\label{THM - M+}
For every $A^+$ as above, $\log p_{A^+}(n) \leq \pi\sqrt{2n|R|/3m}$\;.
\end{theorem}

It is easy to obtain (\ref{eqNathanson}) from the upper bound given by Theorem \ref{THM - M+}. Indeed, we first note that for every $n'$ we have $p_{R^+}(n') \leq (n'+1)^{|R|}$, where $R^{+}=R\setminus \{0\}$.
This follows immediately from the fact that in every partition of $n'$, each of the integers of $R^{+}$ is used at most $n'$ times.
We thus infer that
$$
p_A(n) = \sum_{0 \leq n' \leq n} p_{R^{+}}(n') \cdot p_{A^+}(n-n') \leq (n+1)^{|R|} \sum_{0 \leq n' \leq n} e^{c\sqrt{n-n'}} \leq (n+1)^{|R|+1} e^{c\sqrt{n}} \;,
$$
where $c=\pi\sqrt{2|R|/3m}$. Taking logs from both sides, we obtain (\ref{eqNathanson}).

The proof of Theorem \ref{THM - M+} appears in the next section. At the end of that section we briefly explain why our proof is simpler
than that of \cite{Nat1999}.

\section{Proof of Theorem \ref{THM - M+}.}

For a given fixed integer $m \geq 1$ and $R \subseteq \{0,1,\ldots,m-1\}$, let
$A^{+}$ denote the set of all integers $a \geq m$ with $a \;(\mathrm{mod}\; m) \in R$.
We start with a few observations that extend those used in \cite{Erdos}.
We first note that, for every $0< t < 1$, we have
\begin{equation}\label{eq1}
\sum_{a\in A^{+}} at^a = \sum _{r\in R} \frac{(r+m)t^{r+m}-rt^{2m+r}}{(1-t^m)^2}\;.
\end{equation}

\noindent Indeed, $\sum_{a\in A^{+}} at^a=\sum_{r\in R} \sum _{a\in A_r^{+}}at^a$ where $A_r^{+}$ is the set of all integers $a\geq m$ with $a =r\;(\mathrm{mod}\; m)$ (i.e., $A_r^{+}=\{r+m,r+2m,r+3m,\ldots\}$). Hence, without loss of generality we may assume $|R|=1$. Letting $r\in R$, we have
\begin{align*}
        \sum_{a\in A_r^+} at^a&=t\sum_{a\in A_r^{+}} \frac{d}{dt} t^a=t\cdot\frac{d}{dt}\sum_{a\in A_r^{+}}  t^a =t \cdot \frac{d}{dt} \frac{t^{r+m}}{1-t^m}=\frac{(r+m)t^{r+m}-rt^{2m+r}}{(1-t^m)^2}\;.
\end{align*}
This proves (\ref{eq1}).
We next claim that, if $0 \leq r\leq m-1$ is an integer, then for all $x>0$, we have
\begin{equation}\label{eq2}
    \frac{(r+m)e^{-(r+m)x}-re^{-(2m+r)x}}{(1-e^{-mx})^2}\leq \frac{1}{mx^2}\;.
\end{equation}

\noindent Indeed, since $x>0$, the power series expansion of $e^x$ gives
$$
e^{x/2}-e^{-x/2}=2\sum _{k=0}^{\infty}\frac{1}{(2k+1)!}\left(\frac{x}{2}\right)^{2k+1}=x+x^3\sum_{k=1}^{\infty}\frac{x^{2k-2}}{(2k+1)!\cdot 2^{2k}} > x\;,
$$
implying that
$$
\frac{e^{-x}}{(1-e^{-x})^2}=\frac{1}{(e^{x/2}-e^{-x/2})^2}<1/x^2\;.
$$
We can thus infer that
\begin{align*}
    \frac{(r+m)e^{-(r+m)x}-re^{-(2m+r)x}}{(1-e^{-mx})^2}&={((r+m)e^{-rx}-re^{-(m+r)x}})\frac{e^{-mx}}{(1-e^{-mx})^2}\\
                                                                   &\leq ({(r+m)e^{-rx}-re^{-(m+r)x})}\frac{1}{m^2x^2}\;.
\end{align*}
It remains to check that the expression in parentheses is bounded by $m$. Since the derivative of ${(r+m)e^{-rx}-re^{-(m+r)x}}$ (which is $r(r+m)(e^{-(m+r)x}-e^{-rx})$) is always nonpositive for $x\geq 0$, it is enough to check its value at $x=0$ where it attains the value $m$.
This proves (\ref{eq2}).

We now note that (\ref{eq1}) and (\ref{eq2}) imply that, for every $x>0$,
\begin{equation}\label{eq3}
    \sum_{a\in A^{+}} ae^{-ax} \leq \frac{|R|}{mx^2}\;.
\end{equation}

The final observation we will need is the well-known fact that, for every set of positive integers $S$, we have
\begin{equation}\label{eq4}
    n\cdot p_{S}(n)= \sum_{s\in S\cap [n]} s\sum _{1 \leq k \leq n/s} p_{S}(n-sk)\;,
\end{equation}
where we use $[n]$ for the integers $\{1,\ldots,n\}$. To see this, let $p_{S}(n,s,t)$ and $p_{S}'(n,s,t)$ be the number of partitions of $n$ with summands in $S$ where $s$ appears exactly $t$ times, and at least $t$ times, respectively. Then by double counting,\footnote{The two sides of the first equality count the sum of all integers that appear in all partitions of $n$ using integers from $S$ (there are $p_S(n)$ such partitions). As to the third equality, it follows by observing that each partition of $n$ with exactly $t$ occurrences of $s$ contributes $1$ to $t$ of the summands $p_S'(n,s,t)$, namely $p_S'(n,s,1),p_S'(n,s,2),\ldots,p_S'(n,s,t)$. See Theorem 15.1 in \cite{Nat2000c} for a full detailed proof.} we have
\begin{align*}
    n\cdot p_S(n)&=\sum_{s\in S, t\in \mathbb {N}} s\cdot t \cdot p_S(n,s,t)=\sum _{s\in S\cap [n]}s\sum_{t\in \mathbb {N}}t\cdot p_S(n,s,t)\\
    &=\sum _{s\in S\cap [n]}s\sum_{t\in \mathbb {N}}p_S'(n,s,t)=\sum _{s\in S\cap [n]}s\sum _{1 \leq k\leq n/s}p_S(n-sk)\;.
\end{align*}
This proves (\ref{eq4}).

We are now ready to complete the proof of Theorem \ref{THM - M+}.
We use induction on $n$, with the base case trivially holding. We have
\begin{align*}
    n\cdot p_{A^{+}}(n) &= \sum_{a \in A^+\cap [n]} a\sum _{1 \leq k\leq n/a} p_{A^{+}}(n-ak)\leq \sum_{a\in A^{+}\cap [n]} a\sum _{1 \leq k\leq n/a} e^{c\sqrt{n-ak}}\\
    &\leq e^{c\sqrt{n}}\sum_{a\in A^+\cap [n]} a\sum _{1 \leq k\leq n/a} e^{-\frac{cak}{2\sqrt{n}}}\leq e^{c\sqrt{n}}\sum _{k=1}^\infty \sum_{a\in A^{+}} ae^{-\frac{cak}{2\sqrt{n}}}\\
    &\leq e^{c\sqrt{n}} \sum_{k=1}^{\infty}\frac{4|R|n}{mc^2k^2}= ne^{c\sqrt{n}}\frac{4|R|}{mc^2} \sum_{k=1}^{\infty}\frac{1}{k^2}=n \cdot e^{c\sqrt{n}}\;,
\end{align*}
where the first equality is (\ref{eq4}), the first inequality is by the induction hypothesis, the second inequality uses the elementary fact
$\sqrt{n-rk}\leq \sqrt{n}-\frac{rk}{2\sqrt{n}}$, and in the last inequality we applied (\ref{eq3}) with $x=\frac{ck}{2\sqrt{n}}$.
Dividing both sides by $n$ we obtain the theorem.

\paragraph{Bounding $p_{A^{+}}(n)$ vs. bounding $p_{A}(n)$.}
The reader might be wondering why bounding $p_{A^{+}}(n)$ is so much easier than bounding $p_{A}(n)$. The answer is that the former gives us inequality \eqref{eq2} from which
we obtain the clean inequality \eqref{eq3}. To illustrate the complication that arises when working with $p_{A}(n)$, let us take $A$ to be
the set of odd integers. Then, running the same argument, instead of \eqref{eq2}, one would have liked to use the inequality $\frac{e^{-x}+e^{-3x}}{(1-e^{-2x})^2} \leq \frac{1}{2x^2}$, which is false.
To overcome this, one then needs to use the fact that this inequality is approximately correct for small $x$, which significantly complicates the proof.


\newpage

\noindent{\bf Acknowledgement:}
We would like to thank the referees for their detailed and helpful comments.

\end{document}